\documentclass[10pt,a4paper, leqno]{article}
\usepackage[utf8]{inputenc}
\usepackage{amsmath}
\usepackage{amsfonts}
\usepackage{amssymb}
\usepackage{comment}
\usepackage{hyperref}
\hypersetup{hidelinks
}
\usepackage{pstricks,xy,pst-node}
\xyoption{all}
\usepackage{amsthm}
\usepackage{booktabs} %to produce nicer tables
\usepackage{mathrsfs} %to use \mathscr
\usepackage{mathtools} %to use labelled arrows
\usepackage[disable] %to do notes in the document
{todonotes}
\usepackage[normalem]{ulem} %to strike out text
\usepackage{caption} %to put some space between tables and caption
\usepackage{enumitem} % to change indent of itemize
\usepackage{tikz} %figures
\makeatletter
\let\c@table\c@figure % for (1)
\let\ftype@table\ftype@figure % for (2)
\makeatother

\numberwithin{equation}{section}

\newtheorem{theorem}[equation]{Theorem}
\newtheorem{lemma}[equation]{Lemma}
\newtheorem{proposition}[equation]{Proposition}
\newtheorem{corollary}[equation]{Corollary}

\newtheorem*{theorem*}{Theorem}
\newtheorem*{proposition*}{Proposition}
\theoremstyle{definition}

\newtheorem{remark}[equation]{Remark}
\newtheorem*{ack}{Acknowledgements}
\newtheorem*{notation}{Notation}

\theoremstyle{remark}
\newtheorem{example}[equation]{Example}

\newcommand{\PP}{\mathbb{P}}

\def\cC{{\mathcal C}}

\def\cU{{\mathcal U}}

\def\PP{\mathbf P}

\def\CC{\mathbb{C}}

\def\ZZ{\mathbb{Z}}

\def\cE{{\mathcal E}}

\def\cO{{\mathcal O}}

\def\cQ{{\mathcal Q}}
\def\cL{{\mathcal L}}
\def\ra{\rightarrow}
\def\lra{\longrightarrow}

\def\af1{\mathbf{aff}_1}

\def\KK{\mathrm{Kum}^2(A)}
\def\KN{\mathrm{Kum}^{2,3}(A)}
\def\KU{\mathrm{Kum}^{1}(A)}

\DeclareMathOperator{\rank}{rank}

\DeclareMathOperator{\Sym}{Sym}
\DeclareMathOperator{\Sing}{Sing}
\DeclareMathOperator{\codim}{codim}
\DeclareMathOperator{\Pic}{Pic}

\DeclareMathOperator{\HHH}{H}

\DeclareMathOperator{\Ker}{Ker}
\DeclareMathOperator{\Gr}{Gr}

\DeclareMathOperator{\modulo}{mod}

\DeclareMathOperator{\SU}{SU}
\DeclareMathOperator{\Hilb}{Hilb}

\author{Vladimiro Benedetti\thanks{D\'epartement de math\'ematiques et applications, ENS, CNRS, PSL University, 45 Rue d'Ulm, 75005 Paris, France. ORCID iD: 0000-0001-7113-1639}, Laurent Manivel\thanks{Institut de Math\'ematiques de Toulouse, UMR 5219, Universit\'e de Toulouse, CNRS, UPS IMT, 118 Route de Narbonne, F-31062 Toulouse Cedex 9, France. ORCID iD: 0000-0001-6235-454X}, Fabio Tanturri\thanks{Laboratoire Paul Painlev\'e, UMR CNRS 8524, Universit\'e de Lille, B\^atiment M2, Cit\'e Scientifique, 59655 Villeneuve d'Ascq CEDEX, France. ORCID iD: 0000-0001-5130-9698}}
\title{The geometry of the Coble cubic\linebreak and orbital degeneracy loci}
\begin{document}
\date{}
\maketitle

\begin{abstract}
The Coble cubics were discovered more than a century ago in connection with 
genus two Riemann surfaces and theta functions. They have attracted renewed 
interest ever since. Recently, they were reinterpreted in terms of alternating 
trivectors in nine variables. Exploring this relation further, we show how the 
Hilbert scheme of pairs of points on an abelian surface, and also its Kummer fourfold,
a very remarkable hyper-K\"ahler manifold, 
can very naturally be constructed in this context. Moreover, we explain how this 
perspective allows us to describe the group law of an abelian surface, in a 
strikingly similar way to how the group structure of a plane cubic can be defined 
in terms of its intersection with lines. 
\end{abstract}

\section{Introduction}

The Coble hypersurfaces are very remarkable cubics and quartics in complex projective spaces,
discovered by Coble more than a century ago. They can be characterized as the unique 
hypersurfaces whose singular locus is the Jacobian of a genus two curve embedded in $\PP^8$, 
or the associated Kummer variety of a genus three curve embedded in $\PP^7$, respectively.
\renewcommand{\thefootnote}{\fnsymbol{footnote}} 
\footnotetext{Mathematics Subject Classification: 14K99, 14N05, 14M15, 14J35, 14H60}     
\renewcommand{\thefootnote}{\arabic{footnote}}

The Coble hypersurfaces have been revisited several times. In the eighties, Narasimhan and Ramanan
interpreted them in terms of moduli spaces of vector bundles with fixed determinant on a curve of
genus two or three \cite{NR}. This perspective has been explored by a number of authors, see 
\cite{beauville} and the references therein. 

More recently, the Coble hypersurfaces have been given interpretations coming from Lie theory, 
more precisely from the Kac--Vinberg theory of so-called $\theta$-groups \cite{GSW}. It seems
to be just a coincidence that $\theta$-groups were coined this way at a time where there was
no apparent relation with theta functions, but the fact is that there is a very rich interplay 
between the invariant theory of $\theta$-representations and certain moduli spaces of polarized
abelian varieties. 

From this point of view, genus two curves are naturally related to alternating trivectors,
that is, elements of $\wedge^3V_9$ for $V_9$ a nine-dimensional vector space. Over the 
complex numbers, it was explained in \cite{GSW} how to associate to a general such trivector
an abelian surface $A$ in $\PP(V_9^\vee)$, the projective space of hyperplanes in $V_9$, and 
a cubic hypersurface which is singular exactly along $A$; to be precise, $A$ is only a torsor
over an abelian surface: to make it an abelian surface stricto sensu, one needs to fix an origin. 
%% HERE %%
As further showed in \cite{GS}, the trivector allows one to associate to each point of $A$ a genus two curve $C$ lying on $A$, inducing an isomorphism between $A$ and $A^\vee$. The divisor corresponding to $C$ is a theta divisor on $A$, allowing the identification of $A$ with the Jacobian of $C$; the aforementioned cubic hypersurface has to be the same as the one discovered by Coble. More background material is provided in Section \ref{classicalfacts}; this point of view was further explored over an arbitrary field in \cite{RS}.

In this paper we shall enrich the picture by passing to the dual projective space $\PP(V_9)$,
where our trivector defines a wealth of interesting subvarieties as
orbital degeneracy loci, following the terminology of \cite{BFMT1,BFMT2}. After some preparation in Section \ref{affinemodel}, we investigate their relations with the abelian surface and the Coble cubic in Sections \ref{secODL} and \ref{secKummer}. These subvarietes can a posteriori and more abstractly be interpreted in terms of the branch locus of the natural morphism from the moduli space of semistable rank $3$ vector bundles on $C$ with trivial determinant to $\PP (V_9)$. Already 
well-known were the branch locus itself, a sextic hypersurface which turns out to be the projective dual to the cubic, and a (singular) subvariety $\Sigma$ of its singular locus which can be identified with the triples of 
degree zero line bundles on $C$ whose product is trivial.

A nice
ingredient from the theory of orbital degeneracy loci is that they usually come with 
simple resolutions of singularities, just like the usual degeneracy loci of morphisms 
between vector bundles. We show that our natural resolution of  $\Sigma$ is nothing else
than the Kummer fourfold of $A$ (Theorem \ref{thmHKisom}). Moreover, we observe that, from the 
orbital degeneracy loci point of view, one can define two natural smooth covers of $\Sigma$, 
generically finite of degree three. We identify these covers, one with the Hilbert scheme of
length two subschemes of $A$ (Theorem \ref{isoThilbert}), the other one with the nested Kummer fourfold 
(Theorem \ref{isoK23}). More precisely, the main results can be summarized in the following theorem, where we denoted by $\cU_{r_i}$ the tautological bundle of rank $r_i$ over the flag variety $F(r_1, \dotsc, r_i,\dotsc, V_9)$.
\begin{theorem*}
	Let $\omega\in \wedge^3 V_9$ be a general alternating trivector and let $A$ be its associated abelian surface in $\PP(V_9^\vee)$. If we regard $\omega$ as a general section of the trivial vector bundle $\wedge^3 V_9$ over $\PP(V_9)$, then:
\begin{itemize}
	\item the zero locus of the section induced by $\omega$ of the vector bundle \[\pi^*(\wedge^3 V_9)/(\cU_1\wedge \wedge^2 V_9+ \cU_3\wedge \cU_6\wedge V_9+ \wedge^3 \cU_6)\] over the flag variety $F(1,3,6,V_9)\xrightarrow{\pi} \PP(V_9)$ is isomorphic to the generalized Kummer fourfold $\KK$;

	\item the zero locus of the section induced by $\omega$ of the vector bundle 
	\[\pi^*(\wedge^3 V_9)/(\cU_1\wedge \wedge^2 V_9+ \wedge^2\cU_5\wedge V_9+ \wedge^3 \cU_7)\] over the flag variety $F(1,5,7,V_9)\xrightarrow{\pi} \PP(V_9)$ is isomorphic to the Hilbert scheme $\Hilb^2 (A)$;
	
	\item the zero locus of the section induced by $\omega$ of the vector bundle 
	\[
	\pi^*(\wedge^3 V_9)/(\cU_1\wedge \wedge^2 V_9+ \cU_3\wedge \cU_5\wedge V_9+ \cU_3\wedge \cU_6\wedge \cU_7 + \wedge^3 \cU_6)\] over the flag variety $F(1,3,5,6,7,V_9)\xrightarrow{\pi} \PP(V_9)$ is isomorphic to the nested Kummer fourfold $\KN$.
\end{itemize}
\end{theorem*}

%Finally, we use the same perspective in order to give the following nice description of the group 
%structure on $A$ (defined once an origin $O$ has been fixed):
On the one hand, the above constructions provide to interesting and well-studied objects such as $\KK$ and $\Hilb^2 (A)$ an interpretation as zero loci of sections of suitable vector bundles over some flag varieties. On the other hand, this perspective allows us to give the following nice description of the group structure on $A$ (defined once an origin $O$ has been fixed):
\begin{proposition*}[Proposition \ref{lemmaPQR}]
Let $A \subset \PP(V_9^\vee)$ be the abelian surface associated to a general alternating trivector $\omega \in \wedge^3 V_9$. Then we can fix the origin $O$ in $A$ in such a way that the following holds: for any three general points $P, Q, R\in A \subset \PP(V_9^\vee)$, $P+Q+R=O$ if and only if contracting $\omega$ with any two of the three points yields the same line in $V_9$, i.e., if and only if
	\[
	[\omega(P, Q,\cdot )]=[\omega(P, R,\cdot )]= [\omega(Q, R,\cdot )]\in \PP(V_9).
	\]
\end{proposition*}
This description is, at least 
formally, completely similar to the classical description of the group structure over
a plane cubic, from its intersection with lines. The main difference is that the 
space of ``lines'', rather than the dual projective plane, is now the Kummer fourfold
itself. In terms of the curves $C_P, C_Q, C_R$ associated to $P,Q,R$, the above proposition can be reformulated by saying that $P+Q+R=O$ if and only if $C_P \cup C_Q \cup C_R=H \cap A$ for some hyperplane $H \subset \PP(V_9^\vee)$.

\begin{ack}
	The authors wish to thank Laurent Gruson for stimulating discussions, and Jerzy Weyman for communicating 
	\cite{weymanE8} to them. They moreover thank the referee for the interesting suggestions provided and for pointing out Proposition \ref{lemmaPQRcurves}. \newline
	The third author is supported by the Labex CEMPI (ANR-11-LABX-0007-01).
\end{ack}

\section{Classical facts}
\label{classicalfacts}

\subsection{The Coble cubic}

Let $A$ be an abelian surface, and $\Theta$ a principal polarization. Then $3\Theta$ 
defines an embedding of $A$ inside ${|3\Theta|}^\vee=\PP(V_9^\vee) \simeq \PP^8$, 
where $V_9:= \HHH^0(A,3\Theta)$. The following result is essentially due to Coble \cite{coble}:

\begin{theorem}
\label{coble}
There exists a unique cubic hypersurface $\cC_3$ in $\PP(V_9^\vee)$ which is 
singular along $A$.
\begin{proof}
	See, e.g., \cite[Proposition 3.1]{beauville}.
\end{proof}
\end{theorem}

We will refer to $\cC_3$ as the Coble cubic. 
Note that $A[3]$, the finite group of three-torsion points in $A$, fixes $3\Theta$. 
Therefore  it acts on $\PP(V_9^\vee)$ by leaving $A$ invariant.

\subsection{Moduli of vector bundles on genus two curves}
\label{modvecbun}
Let $C$ be a genus two curve whose Jacobian $JC\cong A$. 
Let $\SU_C(r)$ denote the moduli space of semistable rank $r$ vector bundles on $C$ with trivial determinant. There is a natural morphism from $\SU_C(r)$ to the linear system
%% HERE %%
$|r\Theta|$ (see \cite{Or}):
\begin{itemize}
\item if $r=2$, we get an isomorphism $\SU_C(2)\cong \PP^3$;
\item if $r=3$, we get a finite morphism of degree two 
%% HERE %%
$\SU_C(3)\to |3\Theta|=\PP (V_9)$, branched along a sextic hypersurface $\cC_6\subset \PP (V_9)$.
\end{itemize}

The following result was conjectured by Dolgachev, and proved in \cite{Or} and \cite{Ng}:
\begin{theorem}\label{duality}
The sextic hypersurface $\cC_6$ is the projective dual of the Coble cubic $\cC_3$.
\end{theorem}

\begin{remark}
\label{remsigmainP}
The morphism $\SU_C(3)\to \PP (V_9)$ allow us to describe many subvarieties of $\cC_6\subset \PP (V_9)$ in terms of $\SU_C(3)$. The singular locus of $\cC_6$ is the same as $\Sing(\SU_C(3))$, and can be identified 
with the set of strictly semistable vector bundles on $C$, i.e., vector bundles whose associated graded is $F \oplus \wedge^2(F^\vee)$ for some rank two degree zero vector bundle $F$. Its dimension is five.

Let $A^{(3)}:=\Sym^3 A$, let $\sigma:A^{(3)}\to A$ be the sum morphism, and $\Sigma$ the zero fiber. Then
\[
\Sigma\cong \{E\in \SU_C(3)\mbox{ s.t.\ }E=L_1\oplus L_2\oplus L_3\mbox{ with } L_1, L_2, L_3\in JC\}
\]
has dimension four and can be identified with a codimension one subvariety of $\Sing(\SU_C(3))$.

Finally, the subvariety of $\Sigma$ where two $L_i$'s are isomorphic has dimension two. When all $L_i$'s are isomorphic, we get a zero-dimensional scheme, consisting of 81 points.
\end{remark}

\subsection{Alternating trivectors}
\label{altTriv}

Following \cite{GSW}, one can give another description of the embedding of $A$ in $\PP(V_9^\vee)$, starting from an alternating trivector (or three-form). 

Let $\omega\in \wedge^3 V_9$ be a general alternating trivector. 
Let $H$ denote the hyperplane bundle 
on $\PP(V_9^\vee)$. Then $\wedge^3H$ is a subbundle of the trivial bundle with fiber 
$\wedge^3 V_9$, and the quotient is $\wedge^2 H(1)$. So $\omega$ defines a 
section of $\wedge^2H(1)$ over $\PP(V_9^\vee)$, and the latter can be stratified 
by the rank of this two-form. We denote 
\[
D_{Y^{Sp}_r}:=\{ P\in \PP(V_9^\vee)\mbox{ s.t.\ }\rank(\omega(P,\cdot,\cdot))\leq r\}.
\]
These loci are nothing more than the degeneracy loci (or Pfaffian loci) of the skew-symmetric morphism $H^\vee \ra H(1)$ corresponding to $\omega$. For dimensional reasons, $D_{Y^{Sp}_2}$ is empty, and therefore $D_{Y^{Sp}_4}$ is a smooth surface. By \cite[Theorem 5.5]{GSW}, $D_{Y^{Sp}_4}$ is a torsor, that we denote by $A$, 
over an abelian surface. By \cite[Proposition 5.6]{GSW}, the restriction of the ambient 
polarization is of type $(3,3)$. 
Moreover, the surface $A$ is the singular locus of $D_{Y^{Sp}_6}$, the Pfaffian cubic hypersurface. By Theorem \ref{coble}, this 
hypersurface must be the Coble cubic $\cC_3$. 
Of course all these loci depend on $\omega$, but 
we will omit this dependence in our notation. By varying $\omega$, one gets a locally complete  family of $(3,3)$-polarized 
abelian surfaces \cite{GSW}. The Coble cubic can be also described in terms of an alternating trivector as the fundamental locus of the congruence of lines that the trivector determines (\cite{dpfmr}, see also \cite{tan15} for more on complexes of lines).  

The geometry of the pair $(A,\cC_3)$ was described in more details in \cite{GS}. One remarkable construction is the following: to each point $P \in A$ we associate the kernel of the skew-symmetric matrix $\omega(P,\cdot,\cdot)$, which has (affine) dimension five. The corresponding $\PP^4$ is contained in $\cC_3$; this correspondence induces an isomorphism between $A$ and the Fano variety of $4$-dimensional subspaces of $\cC_3$. Moreover, each such $\PP^4$ cuts $A$ along a theta-divisor \cite[Theorem 3.6]{GS}, inducing the isomorphism between $A$ and $A^\vee$ mentioned in the introduction.

\medskip
On the dual space $\PP(V_9)$, notice that $\omega\in \wedge^3 V_9\cong \HHH^0(\PP(V_9), \wedge^3 \cQ)$, where $\cQ$ denotes the tautological quotient bundle (of rank eight). The fiber of $\wedge^3 \cQ$ is isomorphic to $\wedge^3 \CC^8$. Just as we did on $\PP(V_9^\vee)$ when we defined the Pfaffian loci, we can define subvarieties of 
$\PP (V_9)$ as loci where the trivector that we obtain on $\cQ$ has some special behavior, in the sense that it belongs to some proper $GL_8$-orbit 
(or rather, orbit closure) in $\wedge^3 \CC^8$. This is precisely the idea behind the notion of orbital degeneracy loci introduced in 
\cite{BFMT1,BFMT2}. The next section will be devoted  to the study of the relevant orbits in $\wedge^3 \CC^8$; in the last sections, we will provide geometric interpretations for the corresponding orbital degeneracy loci in $\PP (V_9)$, linking them to the subvarieties described in Remark \ref{remsigmainP}.

\section{The affine model: trivectors in eight variables}

\label{affinemodel}

Our model will be the $GL_8$-representation $\wedge^3 V_8$. This is a classical example
of a representation with a finite number of orbits. The properties we will need in the following regarding this space can be found in \cite{Gurevich,weymanE8}.

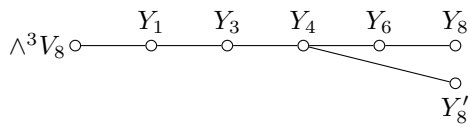
\begin{figure}[h!bt]
 \caption{{\small Inclusions of orbit closures up to codimension eight in $\wedge^3 V_8$}}
  \label{Figure0}
\begin{center}
 \begin{tikzpicture}
    \tikzstyle{every node}=[draw,circle,fill=white,minimum size=4pt,
                            inner sep=0pt]
\path (0,0) node(0) [label=left:$\wedge^3 V_8$] {}
      (1,0) node(1) [label=above:$Y_{1}$] {}
      (2,0) node(3) [label=above:$Y_{3}$] {}
      (3,0) node(4) [label=above:$Y_{4}$] {}
      (4,0) node(6) [label=above:$Y_{6}$] {}
      (5,0) node(8) [label=above:$Y_{8}$] {}
      (5,-0.5) node(88) [label=below:$Y_{8}'$] {};
\draw (0) -- (1);
\draw (1) -- (3);
\draw (3) -- (4);
\draw (4) -- (6);
\draw (6) -- (8);
\draw (4) -- (88);

\end{tikzpicture}
\end{center}
\end{figure}

We denote by $Y_i$ a codimension $i$ orbit closure inside $\wedge^3 V_8$. As we want to construct some orbital degeneracy loci (see Section \ref{ODLbrief}) inside $\PP(V_9)$, whose dimension is eight, we will focus on the varieties $Y_i$ for $i\leq 8$. As it turns out, there is exactly one orbit closure of codimension $i$ for $i=1,3,4,6$ and two distinct orbit closures $Y_8$ and $Y_8'$ of codimension eight.
The inclusion diagram is given in Figure \ref{Figure0}. 

\subsection{Kempf collapsings}
In \cite{weymanE8}, the geometry of these orbit closures has been studied with the 
help of birational Kempf collapsings. These are particular resolutions of singularities given by total spaces 
of homogeneous vector bundles on some auxiliary flag manifolds.

Let $F=G/P$ such a flag manifold, for $P$ a parabolic subgroup of an algebraic group $G$. Then a homogeneous bundle 
on $F$ is of the form $\cE_U=G\times_P U$ for some $P$-module $U$. If $U$ is a $P$-submodule of a $G$-module $V$, then $\cE_U$ is a sub-vector bundle of $\cE_V$, which is the trivial
bundle on $F$ with fiber $V$. In particular, if we denote by $W$ the total space of $\cE_U$, 
this construction induces a proper $G$-equivariant 
map $\pi_U: W\rightarrow V$, called a Kempf collapsing. When $V$ has only finitely 
many $G$-orbits (e.g., the $GL_8$-representation $\wedge^3 V_8$), the image of $\pi_U$ must be some orbit closure $Y$. In many cases, $Y$ being given, we can always find a parabolic $P$ and a $P$-module $U$ such that the image of $\pi_U$
is $Y$.

Partially following \cite{weymanE8}, in Table \ref{kempfTable} we provide a few finite Kempf collapsings for the biggest orbit closures in $V=\wedge^3 V_8$, together with their degrees. We denote by $F(r_1,\ldots,r_i,\ldots,V_8)$ the variety
parametrizing flags of subspaces of $V_8$ of dimensions $r_1,\ldots,r_i$. On this flag
manifold, we denote by $\cU_{r_i}$ the tautological bundle of rank $r_i$. 

	\begin{table}[h!bt]
	\begin{center}
		\caption{Some $d:1$ Kempf collapsings for the orbit closures in $\wedge^8 V_8$ of codimension up to $8$.}
		\begin{small}
		\label{kempfTable}
		\begin{tabular}{cccc} \toprule
		$Y$ &  $F$ & $\mathcal{E}_U$ & $d$ \\
		\midrule
		$Y_1$ & $\Gr(5,V_8)$ & $\wedge^2\cU_5\wedge V_8$ & 1\\
		$Y_3$ & $F(1,4,V_8)$&  $\cU_1\wedge (\wedge^2V_8)+ (\wedge^2\cU_4)\wedge V_8$ & 2\\
$Y_4$ & $F(2,5,V_8)$ & $\wedge^3 \cU_5+\cU_2\wedge \cU_5 \wedge V_8$ & 1 \\
$Y_4$ & $F(4,6,V_8)$ & $\wedge^2 \cU_4 \wedge V_8 + \wedge^3 \cU_6$ & 3 \\
$Y_4$ & $F(2,4,5,6,V_8)$ & $\cU_{2}\wedge(\cU_{4}\wedge V_{8}+\cU_{5}\wedge\cU_{6})+\wedge^3\cU_{5}$ & 3 \\
$Y_6$ & $F(1,3,4,6,7,V_8)$ & $\cU_1\wedge (\cU_4\wedge V_8+\wedge^2 \cU_6)+\cU_3\wedge (\cU_3\wedge \cU_7+\cU_4\wedge \cU_6)$ & 1\\
$Y_8$ & $F(2,5,7,V_8)$& $\cU_2\wedge (\cU_2\wedge V_8 + \cU_5\wedge \cU_7)+\wedge^3 \cU_5$ & 1\\
$Y'_8$ & $\Gr(2,V_8)$ & $\cU_2\wedge (\wedge^2V_8)$ & 1\\
	    \bottomrule
		\end{tabular}
	\end{small}
		\end{center}
	\end{table}

\begin{remark}
\label{nonbiry3}
The Kempf collapsing corresponding to $Y_3$ is finite by a dimension count \cite{weymanE8}. Its degree is at least $2$: indeed, by \cite{Gurevich}, a general element of $Y_3$ is $y_3=v_{123}+v_{456}+v_{147}+v_{268}+v_{358}$ and at least the two flags $\langle v_1 \rangle \subset \langle v_1,v_5,v_6,v_8 \rangle$ and $\langle v_4 \rangle \subset \langle v_2,v_3,v_4,v_8 \rangle$ are in the preimage of $y_3$ in the total space of $\cU_1\wedge (\wedge^2V_8)+ (\wedge^2\cU_4)\wedge V_8$ over $F(1,4,V_8)$. A direct computation in the proof of Proposition \ref{coverbyp3} will show that it is exactly $2$.

In Proposition \ref{propthreetoone} we will show that the second Kempf collapsing for $Y_4$ appearing in Table \ref{kempfTable} is indeed of degree $3$. The third one will appear as the fiber product of the first two at the end of Section \ref{twotriplecovers}.
\end{remark}

\begin{notation}
	For a flag $V_a \subset V_b \subset V_c$, we will write $V_{abc}:=V_a \wedge V_b \wedge V_c$ for short. For instance, for $y \in \wedge^3 V_8$, it turns out that $y \in Y_3$ if and only if $y \in V_{188}+V_{448}$ for some flag $V_1 \subset V_4 \subset V_8$ (see Remark \ref{nonbiry3}), meaning that $y$ can be decomposed into a sum of elements of $V_{188}$ and $V_{448}$.
\end{notation}

A birational Kempf collapsing should be interpreted in the following way: for each point $\omega$
of the open orbit in $Y$, there exists a unique flag $V_\bullet\in F$ such that $\omega$ belongs to 
the fiber of $\cE_U$ over $V_\bullet$. In some sense this yields a normal form for $\omega$. It is mainly from 
this perspective that in the sequel we will make use of Kempf collapsings. 

\medskip 
Let us mention some of the properties of the orbit closures $Y_i \subset \wedge^8 V_8$ that will be useful to us in the sequel.
\begin{itemize}
\item The orbit closures $Y_1, Y_3, Y_4, Y_8$ are normal,
Cohen--Macaulay, and have rational singularities; only $Y_1$ and $Y_4$ are Gorenstein. $Y_6$ and $Y'_8$ are neither normal nor Cohen--Macaulay.
\item Any other orbit of codimension higher than eight is contained in $Y_8 \cup Y'_8$.
\item $Y_1$ is the hypersurface defined by the hyperdeterminant, the unique 
$SL_8$-invariant polynomial of degree $16$ over $\wedge^3 V_8$.
\item $Y_3$ is the singular locus of $Y_1$.
\item The structure sheaf of $Y_4$ admits the following self-dual resolution:
\[
0\to  \det(V_8^\vee)^9 \to \Sym^2 V_8^\vee\otimes \det(V_8^\vee)^5 \to \wedge^4 V_8^\vee\otimes \det(V_8^\vee)^4 \to
\]
\[
\to \Sym^2 V_8 \otimes \det(V_8^\vee)^4 \to \cO_{\wedge^3 V_8}\to \cO_{Y_4}\to 0.
\]
\end{itemize}

\subsection{Two triple covers}
\label{twotriplecovers}
A rather delicate but very interesting point is that there are a priori more
Kempf collapsings than orbit closures. It can happen that some orbit closures have 
several resolutions of singularities by different Kempf collapsings. It can also
happen that a Kempf collapsing is not birational onto its image, either because the
dimension drops or, more scarcely, because it has positive degree. Although the latter
phenomenon cannot happen for the Kempf collapsing of a completely reducible homogeneous 
vector bundle (\cite[Proposition 2 (c)]{kempf}), we already met an instance of it in 
Remark \ref{nonbiry3}. The next example will be essential in the sequel:

\begin{proposition}
\label{propthreetoone}
On the flag manifold $F_2:=F(4,6,V_8)$, consider the homogeneous 
vector bundle $\cE_2=\wedge^2 \cU_4 \wedge V_8 + \wedge^3 \cU_6$. Then the Kempf collapsing 
$\pi_2$ of its total bundle $W_2$ is a generically $3:1$ cover of $Y_4$.

\begin{proof}
Recall from Table \ref{kempfTable} that $Y_4$ admits a resolution of singularities $\pi_1:W_1 \ra Y_4$, where $W_1$ is the total space of the vector bundle $\wedge^3 \cU_5+\cU_2\wedge \cU_5 \wedge V_8$ on $F_1:=F(2,5,V_8)$. Being equivariant, it must be an isomorphism on the open orbit $\cO_4$. We will show that
\begin{enumerate}[label=\roman*.]
	\item if $y \in \pi_2(F_2)$ (or, equivalently, if there exist $U_4 \subset U_6 \subset V_8$ such that $y \in \wedge^2 U_4 \wedge V_8 + \wedge^3 U_6$), then there exists a flag $V_2 \subset V_5$ with $V_2\subset U_4\subset V_5\subset U_6$ such that $y \in V_{555}+V_{258}$. In particular, $\pi_2(F_2)\subset Y_4$;
	\item if $y \in \cO_4$, then there exist exactly three flags $U_4 \subset U_6$ such that $y \in U_{448}+U_{666}$. In particular, $\cO_4 \subset \pi_2(F_2)$.
\end{enumerate}
These two claims will imply that $\pi_2(F_2)=Y_4$ and that $\pi_2$ is a $3:1$ over $\cO_4$.

To prove i., we fix $U_4 \subset U_6$ and we consider the space of parameters for $V_2 \subset V_5$ with $V_2\subset U_4\subset V_5\subset U_6$, i.e., $\Gr(2,U_4)\times
\PP(U_6/U_4)$. Over this space, the point $y$ can be regarded as a section of the trivial vector bundle $\wedge^2 U_4 \wedge V_8 + \wedge^3 U_6$; on a point $V_2 \subset V_5$, $y \in V_2 \wedge U_4 \wedge U_8 + \wedge^3 U_6$ if and only if the induced section $\bar{y}$ on the quotient bundle $V_8/V_6 \otimes \wedge^2 (U_4/V_2)$ vanishes. On the points where $\bar{y}$ vanishes, $y$ induces a section of $V_{258} + \wedge^3 U_6$; it will belong to $V_{258} + \wedge^3 V_5$ if and only the induced section on $\wedge^3 U_6/(V_{258} + \wedge^3 U_5)=(U_6/V_5) \otimes \wedge^2(V_5/U_2)$ vanishes. In other words, the points in $\Gr(2,U_4)\times
\PP(U_6/U_4)$ in the zero locus of the section induced by $y$ of the vector bundle 
\[E=\wedge^2(U_4/V_2)^{\oplus 2}\oplus \wedge^2(U_5/V_2)\]
will give the flags we are looking for. A straightforward computation shows that $c_5(E)=1$, 
so there exists at least one solution to our problem, as claimed.

To prove ii., we start with a point $y \in \cO_4$, for which there exists a unique flag $V_2 \subset V_5$ inside $F_1$ such that $y \in V_{258}+V_{555}$. By i., any flag $U_4\subset U_6$ such that $y \in U_{448}+U_{666}$ has to satisfy $V_2\subset U_4\subset V_5\subset U_6$. This can be reformulated as the condition that contracting $y$ with 
any element $\eta\in U_6^{\perp}$, we get an element of $\wedge^2U_4$. Note that 
since $y\in V_2\wedge V_5\wedge V_8 + \wedge^3 V_5$, such a contraction will
belong to $V_2\wedge V_5$. So it belongs to $\wedge^2U_4$ if and only if its image in 
$V_2\otimes (V_5/V_2)$ is contained in $V_2\otimes (U_4/V_2)$. Dualizing, we need that the 
induced morphism from $U_6^{\perp}\otimes V_2^\vee$ to $V_5/V_2$ has rank two, 
which occurs in codimension two, that is, for a finite number of spaces $U_6$. 
Then $U_4/V_2$ (hence $U_4$) is determined as the image of the previous morphism. 
We conclude that the number of solutions to our problem can be computed 
on $\PP(V_8/V_5)$, by the Thom--Porteous formula, as the class
\[
c_2(U_6^{\perp}\otimes V_2^\vee-V_5/V_2) =c_2(U_6^{\perp}\otimes V_2^\vee) =3. 
\]
Being equivariant, $\pi_2$ must be finite and \'etale over the open orbit $\cO_4$ in $Y_4$, hence there exist exactly three flags over any point in the orbit. \qedhere
\end{proof}
\end{proposition}

\begin{example}
Let us denote again by $W_1$ the (total space of the) vector bundle $\wedge^3 \cU_5 + \cU_2 \wedge \cU_5 \wedge \cU_8$ over $F(2,5,V_8)$, which yields a Kempf collapsing resolving $Y_4$.
Consider the point $y_4=v_{456}+v_{147}+v_{257}+v_{268}+v_{358}+v_{367}$ of $\cO_4$, see \cite{Gurevich}.
Its unique preimage in $W_1$ is given by the flag 
\[
V_2=\langle v_7,v_8 \rangle \quad \subset \quad V_5=\langle v_4,v_5,v_6,v_7,v_8 \rangle.
\]
Its three preimages in $W_2$ are given by  
\[\begin{array}{rcl}
U_4=\langle v_5,v_6,v_7,v_8\rangle & \subset & U_6=\langle v_1,v_4,v_5,v_6,v_7,v_8\rangle , \\
U_4=\langle v_4,v_5-v_6,v_7,v_8\rangle & \subset & 
U_6=\langle v_2+v_3,v_4,v_5,v_6,v_7,v_8\rangle , \\
U_4=\langle v_4,v_5+v_6,v_7,v_8\rangle & \subset &
U_6=\langle v_2-v_3,v_4,v_5,v_6,v_7,v_8\rangle .
\end{array}\]
\end{example}

Let $y$ be a point in $Y_4$, and consider a flag $V_2\subset V_4\subset V_5\subset V_6$ 
defining points of $W_1$ and $W_2$ over $y$. This means that $y$ belongs to the intersection 
of $\wedge^3V_5+V_2\wedge V_5 \wedge V_8$ with $\wedge^3V_6+\wedge^2V_4 \wedge V_8$, that is, to 
$V_{248}+V_{256}+V_{555}$.  
This suggests to consider, on the flag manifold $F(2,4,5,6,V_8)$, the vector bundle 
$\cU_{248}+\cU_{256}+\cU_{555}$, and to denote its total space by $W_{12}$. Note that 
the latter bundle has rank $27$ on a $25$-dimensional flag manifold, so that $W_{12}$ has 
dimension $52$, as expected.

We get the following diagram of morphisms:
\begin{equation}
\label{affinediag}
\xymatrix @C=2pc @R=0.8pc{
& & \rule{1pt}{0pt}W_{12} \ar[dl]_{{3:1}} \ar[d] \ar[dr]^{\mbox{\small birat}}  \\
& W_1 \ar[dl]  \ar[dr]_{{\mbox{\small desing}}} & F(2,4,5,6,V_8) & W_2 \ar[dl]^{{3:1}}  \ar[dr] \\
F(2,5,V_8) & & Y_4 \ar@{^{(}->}[d] & & F(4,6,V_8) \\
& & \wedge^3 V_8 & &
} 
\end{equation}
where $W_{12} \ra W_1$ is generically 3:1.

\section{The degeneracy loci}
\label{secODL}

\subsection{Orbital degeneracy loci: generalities}

\label{ODLbrief}

Let us briefly recall how orbital degeneracy loci are constructed \cite{BFMT1, BFMT2}. 
One starts with a model, that we will choose to be a representation $V$ of some 
algebraic group $G$. Inside $V$, we consider a closed $G$-stable subvariety $Y$, usually the closure of a $G$-orbit. For any $G$-principal bundle $\cE$ over a variety $X$, we can consider its associated vector bundle $\cE_V$ on $X$. By construction, each fiber of this bundle gets 
identified with $V$, not canonically, but the ambiguity only comes from the action of $G$.
This allows us to define, for any global section $s$ of $\cE_V$, the orbital degeneracy locus
\[
D_Y(s)=\{x\in X \mbox{ s.t.\ } s(x)\in Y\subset V\simeq \cE_{V,x}\}.
\]
In the algebraic context, there is a natural scheme structure induced on $D_Y(s)$ that we will not consider.
In the usual situation where $s$ is general in a finite-dimensional
space of global sections that generates $\cE_V$ everywhere, the orbital degeneracy loci
are well-behaved, in the sense that their properties faithfully reflect the properties of $Y$. 
In particular, $\Sing(D_Y(s))=D_{\Sing (Y)}(s)$ and
\[
\codim_X D_Y(s) =\codim_V Y, \quad \codim_{D_Y(s)} \Sing(D_Y(s))= \codim_{Y} \Sing(Y).
\]

A remarkable feature of an orbital degeneracy locus $D_Y(s)$ associated to a subvariety $Y$ admitting a Kempf collapsing is that it is easy to relativize such a collapsing to get a surjective map $\mathscr{Z} \ra D_Y(s)$. $\mathscr{Z}$ turns out to be the zero locus of an induced section of a vector bundle on a manifold, both determined by the collapsing. Moreover, if the Kempf collapsing is finite, such a map will be finite as well. We refer to \cite{BFMT1, BFMT2} for more details.

\subsection{Loci associated to an alternate trivector}
An example of orbital degeneracy locus is the abelian surface $A=D_{Y^{Sp}_4}$ constructed in Section \ref{altTriv}. Recall that $A$ is in fact only a torsor over an abelian surface, as it has no fixed origin; in order to simplify the terminology, from now on we will abusively call it the abelian surface $A$. Here the model $V$ is the space of alternating bivectors $\wedge^2 V_8$, on which the group $GL_8$ acts. The stable closed subvarieties are the loci $Y^{Sp}_r$ where the rank is bounded above by $r$. The trivector  $\omega\in \wedge^3 V_9$ defines a section of $\wedge^2 H(1)$ over $\PP(V_9^\vee)$, where $H$ denotes the tautological bundle 
on $\PP(V_9^\vee)$, and the corresponding orbital degeneracy loci are the Pfaffian loci $D_{Y^{Sp}_r}\subset \PP(V_9^\vee)$.

On the dual projective space $\PP(V_9)$, 
the trivector  $\omega\in \wedge^3 V_9$ can be seen as a section of $\wedge^3 \cQ$, $\cQ$ being the tautological quotient bundle on $\PP(V_9)$. For $Y_i$ the orbit closures inside $\wedge^3 V_8$ introduced in Section \ref{affinemodel}, the associated degeneracy loci are
\[
D_{Y_i}:=\{ [V_1]\in \PP(V_9) \mbox{ s.t.\ }  \omega \,(\modulo V_1)\in Y_i\subset \wedge^3(V_9/V_1)\},
\]
where we omit, for simplicity, the dependence on $\omega$. We will often write $\omega / V_1$ instead of $\omega \,(\modulo V_1)$ where no confusion can arise.

For a generic $\omega$, $D_{Y_i}$ has codimension $i$ inside $\PP(V_9)$ and $\Sing(D_{Y_i}) =D_{\Sing(Y_i)}$ (see Figure \ref{Figure1} for the inclusion graph). For example, $D_{Y_1}$ is a sextic hypersurface inside $\PP(V_9)$, singular along $D_{Y_3}$, which is five-dimensional.

\begin{figure}[h!bt]
	\caption{{\small Inclusions of degeneraci locy inside $\PP(V_9)$}}
	\label{Figure1}
	\begin{center}
		\begin{tikzpicture}[scale=1.3]
		\tikzstyle{every node}=[draw,circle,fill=white,minimum size=4pt,
		inner sep=0pt]
		\path (0,0) node(0) [label=left:$\PP(V_9)$] {}
		(1,0) node(1) [label=above:$D_{Y_{1}}$] {}
		(2,0) node(3) [label=above:$D_{Y_{3}}$] {}
		(3,0) node(4) [label=above:$D_{Y_{4}}$] {}
		(4,0) node(6) [label=above:$D_{Y_{6}}$] {}
		(5,0) node(8) [label=above:$D_{Y_{8}}$] {}
		(5,-0.5) node(88) [label=below:$D_{Y_{8}'}$] {};
		\draw (0) -- (1);
		\draw (1) -- (3);
		\draw (3) -- (4);
		\draw (4) -- (6);
		\draw (6) -- (8);
		\draw (4) -- (88);
		
		\end{tikzpicture}
	\end{center}
\end{figure}

\begin{proposition}
$D_{Y_1}$ is the Coble sextic $\cC_6$.
\begin{proof}
Let us prove that $D_{Y_1}$ is the dual hypersurface to the Coble cubic 
$\cC_3=D_{Y^{Sp}_6}$. The conclusion will then follow from Theorem \ref{duality}.

A general point of the cubic hypersurface $D_{Y^{Sp}_6}$ is a hyperplane $H$
such that we can decompose $\omega = e\wedge\theta +\sigma$ with $e\notin H$, 
$\sigma\in\wedge^3H$, and $\theta\in\wedge^2H$
is degenerate, i.e.\ $\theta^4:=\theta \wedge \theta \wedge \theta \wedge \theta=0$. At a smooth point of this hypersurface, $\theta$ has rank 
six. Let us analyze how $H$ can be deformed inside $D_{Y^{Sp}_6}$. 

We choose a basis $e_1,\ldots ,e_9$ of $V_9$ such that $H$ is generated by $e_1, \ldots , e_8$,
and $e_9=e$. On a neighborhood of $H$ in $\PP(V_9^\vee)$, a hyperplane $H_z$ has a basis $f_1=e_1+z_1e_9, \ldots , f_8=e_8+z_8e_9$. In the basis $f_1,\ldots , f_8, e_9$ of $V_9$, our 
$\omega$ decomposes as $e_9\wedge (\theta+z\lrcorner\sigma)(f)+\sigma(f)$, where 
$z=z_1e^*_1+\cdots +z_8e_8^*$ and 
$\theta(f), \sigma(f)$ are obtained by expressing $\theta, \sigma$ in the basis $e_1, \ldots , e_8$, and 
replacing formally each $e_i$ by $f_i$. So $H_z$ remains inside $D_{Y^{Sp}_6}$
if and only if $(\theta+z\lrcorner\sigma)^4=0$. In particular, the tangent hyperplane
to $D_{Y^{Sp}_6}$ at $H$ is given by the condition that 
$\theta^3\wedge (z\lrcorner\sigma)=0$. 

Suppose our basis has been chosen so that $\theta = e_{12}+e_{34}+e_{56}$, and decompose further 
our tensor as 
\[\omega =e_9\wedge (e_{12}+e_{34}+e_{56})+e_7\wedge\sigma_7+e_8\wedge\sigma_8+e_7\wedge e_8\wedge \sigma_{78},\]
where $\sigma_{7}, \sigma_{8}, \sigma_{78}$ only involve $V_6 : =\langle e_1, \ldots , e_6
\rangle$. The condition 
$\theta^3\wedge (z\lrcorner\sigma)=0$ simply becomes $z\lrcorner\sigma_{78}=0$, or 
equivalently, that  $\sigma_{78}$ belongs to $H_z$. In other words, the tangent hyperplane 
to $D_{Y^{Sp}_6}$ at $[H]$ is the orthogonal hyperplane to the vector $\sigma_{78}$. 

Finally, we claim that $[\sigma_{78}]$ belongs to $D_{Y_1}$. Indeed, when we 
mod out by $\sigma_{78}$, we get 
\[\bar\omega =e_9\wedge (e_{12}+e_{34}+e_{56})+e_7\wedge\sigma_7+e_8\wedge\sigma_8,\]
now in $\wedge^3(V_9/\langle \sigma_{78}\rangle)$ (with some abuse of notation), 
and the factors of $e_9, e_7, e_8$ 
now live in $\wedge^2(V_6/\langle \sigma_{78}\rangle)$. If  $\bar{V}_5 := V_6/\langle \sigma_{78}\rangle$ and $\bar{V}_8 := V_9/\langle \sigma_{78}\rangle$, we conclude that 
$\bar\omega$ belongs to $\wedge^2\bar{V}_5\wedge \bar{V}_8$. But this is precisely 
the condition that defines $D_{Y_1}$. We have thus proved that the dual of 
$\cC_3=D_{Y^{Sp}_6}$ is contained in $D_{Y_1}$.

\smallskip
Conversely, let us consider a general point $[V_1]$ of $D_{Y_1}$, generated
by $e_1$. From the Kempf
collapsing resolving $Y_1$, we see that this means that there must exist a unique $V_6$, 
with $V_1\subset V_6\subset V_9$, such that $\omega$ belongs to $V_{199}+V_{669}$. 
This is equivalent to the fact that the contraction by $\omega$ 
sends $\wedge^2V_6^\perp$ to $V_1$. 

Let us consider a basis $e_1, \ldots , e_6$ of $V_6$.
Note that  $(V_{199}+V_{669})/V_{669}$ is isomorphic to $V_1\otimes \wedge^2(V_9/V_6)$. 
Since $V_9/V_6$ is three-dimensional, every bivector in $\wedge^2(V_9/V_6)$ is decomposable. 
This allows us to complete our basis of $V_6$ in a basis of $V_9$ with three vectors
$e_7, e_8, e_9$ in $V_9$ such that 
\[\omega = e_{178}+\phi_7\wedge e_7+ \phi_8\wedge e_8+ \phi_9\wedge e_9+\psi,\]
where $\phi_7, \phi_8, \phi_9$ belong to $\wedge^2V_6$ and $\psi$ to to $\wedge^3V_6$. 
We claim that the tangent space to $D_{Y_1}$ at $[V_1]$ is the hyperplane defined 
by the linear form $e_9^*$ from the dual basis. Indeed, we describe points 
$[U_1]$ in $D_{Y_1}$ locally 
around $[V_1]$ by moving $V_6$ to spaces $U_6$ such that the contraction by $\omega$ 
from $\wedge^2U_6^\perp$ to $V_9$ keeps rank one, and defining $U_1$ as the image. 
Locally around $V_6$, such a space $U_6$ is defined by linear forms $f_i^*=e_i^*+t_i$, 
for $i=7, 8, 9$, where $t_i$ is a linear combination of $e_1^*, \ldots , e_6^*$. 
Clearly the contraction $\omega (f_7^*, f_8^*, \cdot)$ is a non-zero vector, which
must therefore generate $U_1$, and its coefficient on $e_9$ is $\phi_9(t_7,t_8)$,
which has order two. Modding out order two deformations, $U_1$ is thus contained in $\langle e_1,\dotsc,e_8 \rangle$, which implies the claim. 

In order to conclude the proof, we just need to check that this tangent hyperplane
belongs to $\cC_3=D_{Y^{Sp}_6}$, or equivalently, that the contraction $\omega (e_9^*,\cdot,
\cdot)$ has rank at most six. But that is clear, since this contraction is $\phi_9$,
an element of $\wedge^2V_6$. 
\end{proof}
\end{proposition}

\begin{corollary}
$D_{Y_3}$ is the singular locus of the sextic $\cC_6$. 
\end{corollary}

\begin{proposition}
\label{propbirA}
There exists a natural birational map $D_{Y_6}\dashrightarrow A$. 

\begin{proof}
Let $[e_0]$ be a point in $D_{Y_6}$. By definition, this means that 
we can decompose $\omega$ with respect to a decomposition $V_9=\CC e_0\oplus V_8$ as
$\omega = e_0\wedge\alpha + \eta$, where $\eta$ belongs to $Y_6\subset \wedge^3V_8$. 
By Table \ref{kempfTable}, for $[e_0]$ outside $D_{Y_8}$, this implies that there exists a unique flag  $V_1\subset 
V_4\subset V_7\subset V_8$ such that $\eta$ belongs to $V_{148}+\wedge^3V_7$. 
Let $W_8=\CC e_0\oplus V_7$. Let us also choose a generator $e_1$ of $V_1$, and 
some $e_8\notin W_8$. We can then write $\omega$ as 
\[\omega = e_0\wedge u\wedge e_8+e_1\wedge v\wedge e_8+\xi,\]
for some $u \in V_7$, $v \in V_4$, $\xi \in \wedge^3W_8$. Since the two-form $e_0\wedge u+e_1\wedge v$ has 
rank at most four, this implies that $W_8$ is a point of $A$. We have thus a rational map sending $[e_0]$ to $V_7([e_0])\oplus\CC e_0$.

\smallskip
Conversely, we claim that for a general $P\in A$, there exists a unique line $l_P\subset P$ 
such that $(\omega / l_P)\in Y_6 \subset \wedge^3(V_9/l_P)$ and $V_7((\omega / l_P))=P / l_P$.
By hypothesis there exist $\sigma\in \wedge^2 P$ of rank four and $\sigma'\in \wedge^3 P$, $v_8\not\in P$ such that $\omega = v_8\wedge \sigma + \sigma'$.

We want to show that there exists in general exactly one flag $V_1 \subset V_3 \subset V_4 \subset V_6 \subset V_7 \subset (V_9/l_P)$ such that $\omega/l_P \in V_{148} + V_{166} + V_{337} + V_{346}$.
The space $V_7=P / l_P$ is determined by $l_P$. Denote by $W$ the four-dimensional space defined by $\sigma$; then, $l_P \subset W$ and therefore $ \PP (W)$ is the right parameter space for $l_P$. Indeed, the contraction of $\omega/l_P$ by any element of $(P/l_P)^\perp$ has to belong to $V_1 \wedge V_4$, hence it has rank at most two, but this means that also the rank of $\sigma / l_P$ can be at most two. Moreover, the two-dimensional space $ l_P^{\perp_\sigma} / l_P$ defined by $\sigma / l_P$ must contain $V_1$ and be contained inside $V_4$. Therefore $V_1$ is parametrized by $\PP (l_P^{\perp_\sigma} / l_P)$, while $V_4 /  l_P^{\perp_\sigma} \subset V_6/  l_P^{\perp_\sigma} \subset P /  l_P^{\perp_\sigma}$ is parametrized by $F (2,4,P /  l_P^{\perp_\sigma})$ and $V_3 / V_1 \subset V_4 / V_1$ is parametrized by $\Gr (2,V_4 / V_1)$.

Inside this fourteen-dimensional parameter space, $(v_8 \wedge \sigma)/l_P$ belongs by construction to $V_{148}$. If we interpret $\sigma' / l_P$ as a section of the bundle $\wedge^3  (P /  l_P)$, the flags we are looking for are defined by the vanishing of the section that $\sigma' / l_P$ induces on the quotient bundle $\wedge^3 V_7/(V_{147}+V_{337}+V_{346}+V_{166})$. To determine the number of such flags, we must compute the top Chern class of this quotient 
bundle. In order to do this, we first filter our bundle  by homogeneous subbundles such that the associated graded 
bundle is completely reducible. Explicitely, the associated bundle we get is 
\begin{multline*}
(V_7 / V_6) \otimes \big(V_4 \otimes (V_6 / V_4) \oplus \wedge^2(V_6 / V_4) \oplus (V_4 / V_3)\otimes (V_3 / V_1)  \big) \oplus
\\
\oplus \wedge^2 (V_6 / V_4)\otimes (V_4 / V_1).
\end{multline*} 
A computation with \cite{Macaulay2} shows that the top Chern class of the latter bundle is $1$, hence also of the original one, and the claim follows.

\smallskip
We have thus defined two rational maps $D_{Y_6}\dashrightarrow A$ and $A\dashrightarrow 
D_{Y_6}$, inverse one to the other. This implies the claim.
\end{proof}
\end{proposition}

We will now focus on the four-dimensional orbital degeneracy locus $D_{Y_4}\subset \PP(V_9)$. A priori, it contains $D_{Y_6}$ (dimension two), $D_{Y_8}$ and $D_{Y_{8}'}$ (dimension 0). Following \cite{BFMT1, BFMT2}, we can relativize the three Kempf collapsings of Table \ref{kempfTable} and Section \ref{twotriplecovers}. We will denote by $Z$ the desingularization of $D_{Y_4}$ which is a zero locus inside the flag bundle $F(2,5,\cQ)\cong F(1,3,6,V_9)$ and which relativizes the first Kempf collapsing; more precisely, it is the zero locus of a section of the vector bundle
\[
\wedge^3 V_9/(\cU_{199}+\cU_{369}+\cU_{666})
\]
%\[ (\cU_6/\cU_1)\otimes \wedge^2(V_9 / \cU_6) + \wedge^3(V_9/\cU_6) + \wedge^2(\cU_6/\cU_3) \otimes V_9/\cU_6 \]
induced by $\omega$. Similarly, the relativization of the second one yields a generically $3:1$ morphism to $D_{Y_4}$ from a variety $T$ which is a zero locus of a section of the vector bundle 
\[
\wedge^3 V_9/(\cU_{199}+\cU_{559}+\cU_{777})
\] 
inside $F(4,6,\cQ)\cong F(1,5,7,V_9)$. Finally, we have also an associated zero locus $X$ of a section of the vector bundle
\[
\wedge^3 V_9/(\cU_{199}+\cU_{359}+\cU_{367}+\cU_{666})
\]
inside the flag bundle $F(2,4,5,6,\cQ)\cong F(1,3,5,6,7,V_9)$ which relativizes the fiber product $W_1\times _{Y_4}W_2$, see \eqref{affinediag}. The situation is described in \eqref{reldiag}.

\begin{equation}
\label{reldiag}
\xymatrix @C=2pc @R=0.8pc{
& & \rule{1pt}{0pt}X \ar[dl]_{{3:1}} \ar@{^{(}->}[d] \ar[dr]^{\mbox{\small birat}}  \\
& Z \ar@{_{(}->}[dl] \ar[dr]_{\mbox{\small desing}} & F(2,4,5,6,\cQ) & T \ar[dl]^{{3:1}}  \ar@{^{(}->}[dr] \\
F(2,5,\cQ) & & D_{Y_4} \ar@{^{(}->}[d] & & F(4,6,\cQ) \\
& & \PP(V_9) & &
} 
\end{equation}

\begin{remark}
The varieties $Z,T,X$ depend on $\omega$. Once again, for the sake of lightness of notation, we will omit to write this dependence explicitly.
\end{remark}

By using the description of the model $Y_4$ given in the previous section, the following 
facts can be checked. 

\begin{proposition} 
\label{propmanyfacts}
Let $\omega\in\wedge^3V_9$ be a generic trivector. Then:
\begin{enumerate}
\item $D_{Y_8}=\{ p_1,\cdots,p_{81} \}$ consists in  $81$ reduced points, while $D_{Y_{8}'}=\emptyset$.
\item The surface $D_{Y_6}$ is smooth outside $D_{Y_8}$.  
\item $D_{Y_4}$ is a normal, linearly non-degenerate fourfold with $h^i(\cO_{D_{Y_4}})=1$ for $i=0,2,4$, and $0$ otherwise.
\item $Z$ has trivial canonical bundle, $\chi(\cO_{Z})=3$, $\chi(\Omega^1_{Z})=-6$ and $\chi(\Omega^2_{Z})=90$.
\item $T$ also has trivial canonical bundle, but $\chi(\cO_T)=\chi(\Omega^1_T)=0$.
\end{enumerate}

\begin{proof} Statements $1.,3.,4.,5.$ can be proved by using the desingularization of the loci given by the Kempf collapsings of the affine model, and by computing the corresponding Chern classes with \cite{Macaulay2}. Statement $2.$ is a consequence of the fact that $\Sing(Y_6)\subset 
Y_8\cup Y'_8$. 
\end{proof}
\end{proposition}

\begin{corollary}\label{hk}
$Z$ is a hyper-K\"ahler fourfold.

\begin{proof} This follows from $4.$ and the Beauville--Bogomolov decomposition. \end{proof}
\end{corollary}

\section{The Kummer geometry of a trivector}
\label{secKummer}

This section relates the loci we have constructed with the geometry of the abelian surface $A=D_{Y^{Sp}_4}$. 
More precisely, in the four theorems of this section we will identify the varieties $D_{Y_4},Z,T,X$ with four ``classical'' fourfolds that can be constructed from $A$.

\subsection{The Kummer fourfold}
We will denote by $\KK$ the generalized Kummer hyper-K\"ahler fourfold associated to $A$. It is a subvariety of the Hilbert scheme $\Hilb^3(A)$ of three points over $A$. Recall that there is a well-defined natural morphism $HC:\Hilb^3 (A) \to A^{(3)}$, called the Hilbert--Chow morphism.
Composing with the sum map $A^{(3)}\ra A$, we get a morphism  $\Hilb^3 (A)\ra A$ whose
fibers are all copies of $\KK$ (note that this ensures that the Kummer fourfold is not 
affected by the choice of the origin in our torsor $A$). The Kummer fourfold is a resolution
of the singularities of the fiber $\Sigma\subset  A^{(3)}$ of the sum morphism.

\begin{theorem}
\label{thmHKisom}
$Z$ is isomorphic to the generalized Kummer fourfold $\KK$.

\begin{proof}
We will construct a finite flat morphism $Z'\ra Z$ of degree three. It will 
induce a morphism from $Z$ to $\Hilb^3(A)$. Since $Z$ is hyper-K\"ahler, 
the composition $Z\ra A$ with the sum map must be constant, otherwise we would get 
non-trivial one-forms on $Z$. Therefore 
our morphism factorizes through $\KK$. Finally, it will turn out to be birational.  
Since $Z$ is a minimal model, such a birational morphism must 
be an isomorphism. 

Recall that $Z$ is embedded inside $F(1,3,6,V_9)$. Denote by $F_Z$ the restriction to $Z$
of the $\PP^2$-bundle defined by the natural projection $F(1,3,6,8,V_9)\ra F(1,3,6,V_9)$.
By definition, for any flag $V_1\subset V_3\subset V_6\subset V_9$ in $Z$, $\omega$ belongs
to $V_{199}+V_{369}+V_{666}$. Modding out the latter bundle by $V_{199}+V_{339}+V_{666}$, we 
get the vector bundle $A_2\otimes B_3\otimes C_3$, where $A_2=V_3/V_1$, $B_3=V_6/V_3$ and 
$C_3=V_9/V_6$. Since $F_Z=\PP_Z(V_9/V_6)^\vee$, the class $\bar{\omega}$ of $\omega$ in $A_2\otimes B_3\otimes C_3$ defines a morphism 
\[\hat{\omega} : A_2^\vee (-1)\lra B_3\] 
of vector bundles over $F_Z$. We define $Z'$ as the first degeneracy locus $D_1(\hat{\omega})$. We will prove that $Z'\ra Z$ is finite of degree three, i.e., the fibers always have expected codimension. As they are determinantal, this implies that they are Cohen--Macaulay, hence the projection $Z'\ra Z$ is flat.
It will further induce a morphism to $\Hilb^3( A)$ because the points in the fibers 
of $Z'\ra Z$ are defined by hyperplanes $V_8$ that must belong to $A$. Indeed, 
since $V_8$ contains $V_6$, when we mod out by $\wedge^3V_8$ the class of $\omega$ 
belongs to $(V_{18}+V_{36})\otimes (V_9/V_8)$. When $V_8$ defines a point of 
$D_1(\hat{\omega})$, the term from
$V_{36}$ has rank two modulo $V_1$, so that the image of $\omega$ in $V_{18}+V_{36}$
has rank at most four, which is exactly the condition for $V_8$ to belong to $A$. 

There remains to check that the projection $Z'\ra Z$ is indeed finite of degree three. 
First note that on $\PP (C_3^\vee)$, the morphism $\hat{\omega}: A_2^\vee (-1)\lra B_3$
is expected to drop rank in codimension two, hence on a finite scheme of length $c_2(B_3-A_2^\vee (-1))=s_2(A_2^\vee (-1))=3$.

Let us prove that $D_1(\hat{\omega})$ cannot be positive dimensional.
By what has been said before, $D_1(\hat{\omega})$ is contained in $A$. Moreover, it is defined 
as a subscheme of a projective plane by three quadrics, the $2\times 2$-minors of the matrix $\hat{\omega}$. If it is not the whole plane, this immediately implies that $D_1(\hat{\omega})$
is contained in a conic. In any case, if $D_1(\hat{\omega})$ has positive dimension, it must contain a rational curve. Since $A$ does not contain any rational curve, we get a contradiction. 

Finally, we have to prove that the morphism is birational; for this sake, we provide an explicit description of the image of the morphism on a general point and show that it is generically injective. Let $[e_0] \in D_{Y_4}$ be a general point. Let $V_9=\CC e_0 \oplus V_8$. Then $\omega=e_0 \wedge \sigma + \sigma'$, where $\sigma \in \wedge^2 V_8$ and
\[
\sigma'=v_{456}+v_{147}+v_{257}+v_{268}+v_{358}+v_{367}
\]
in a suitable basis of $V_8$. One readily checks that the hyperplanes $[v_1^*]^\perp, [v_2^*+v_3^*]^\perp, [v_2^*-v_3^*]^\perp$ belong to $A$ and contain $V_6$, hence they must be
the three points in $A$ which correspond to $[e_0]$ via the morphism. If we contract any two
linear forms among $v_1^*, v_2^*+v_3^*, v_2^*-v_3^*$ with $\sigma'$ we get zero, so the same contraction with $\omega$ yields a multiple of $e_0$ (non-zero, since $\sigma$ is general). 
This means that we can generically recover $[e_0]$ from its image in $\Hilb^3(A)$.
\end{proof}
\end{theorem}

\begin{theorem}\label{sigma}
$D_{Y_4}$ is projectively equivalent to $\Sigma$.  

\begin{proof} We want to compare the vertical projections in the diagram
\begin{equation*}
\label{projections}
\xymatrix @C=2pc @R=0.8pc{
& Z\ar[d]_p & \simeq & \KK\ar[d]^q & \\
\PP (V_9) & D_{Y_4}\ar[l] & & \Sigma\ar[r] & |3\Theta|
} 
\end{equation*}
Recall that since $\KK$ has no holomorphic one-forms, its Picard group and its N\'eron-Severi group are the same. Moreover, by \cite[Proposition 8]{beauville83},
the N\'eron-Severi group of the Kummer fourfold is 
\[NS(\KK)=\Pic(\KK) =\iota (NS(A))\oplus \ZZ E,\]
where the map $\iota$ is injective and $E$ is the exceptional divisor of the 
projection from $\KK$ to 
$\Sigma$. We will denote by $L^{[2]}$ the image of $L\in NS(A)$ by $\iota$.
For $\omega$ generic the abelian surface $A$ is generic, so $NS(A)=\ZZ \Theta$.
The projection $q$ is then defined by the full linear system $|\Theta^{[2]}|$. 

We will show below that $p$ also contracts $E$, so that the pull-back
$\cL$ of the dual tautological line bundle on $\PP (V_9)$ must be of the form 
$L^{[2]}$ for some $L\in \Pic(A)$, hence $\cL=\ell\Theta^{[2]}$. By \cite[Lemma 5.2]{britze}, 
\[\chi (\KK,\ell\Theta^{[2]})=3\binom{\ell+2}{2}.\]
A computation
with \cite{Macaulay2} yields that $\chi (Z,\cL)=9$, hence $\ell=1$. Since $D_{Y_4}$ is linearly
non-degenerate by Proposition \ref{propmanyfacts}, $p$ must be
defined by the full linear system $|\cL|\simeq\PP^8$. So $p$ and $q$ are
the same maps, and the conclusion follows.

It remains to show that $p$ contracts $E$. Recall from the proof of Theorem \ref{thmHKisom} how we constructed an isomorphism from 
$Z$ to $\KK$: for any $V_1$ in $D_{Y_4}$ and any flag $V_1\subset V_3\subset V_6\subset V_9$
such that $\omega$ belongs to $V_{199}+V_{369}+V_{666}$ (hence defining a point of $Z$ above $V_1$), 
we deduced an element $\bar\omega$ of $A_2\otimes B_3\otimes C_3$, where $A_2=V_3/V_1$, $B_3=V_6/V_3$ and 
$C_3=V_9/V_6$. Then we proved that the first degeneracy locus of the induced morphism $A_2^\vee(-1)
\ra B_3$ over $\PP(C_3^\vee)$ defines a length three subscheme of $A$. 

It will sufficient to show that the preimage $p^{-1} (D_{Y_6})$ is a three-dimensional subscheme of $E$. Since $E$ is irreducible, the two are in facts equal and their image through $p$ is $2$-dimensional. If $V_1$ is a general point of $D_{Y_6}$, then we can write $\omega / V_1$ as $v_{456}+v_{147}+v_{257}+v_{268}+v_{358}$ for a suitable choice of a basis of $V_9/V_1$. This determines the unique flag
\begin{multline}
\label{uniqueflag}
\langle v_8 \rangle \subset \langle v_4, v_7, v_8 \rangle \subset
\langle v_4, v_5, v_7, v_8 \rangle \subset  \langle v_2, v_4, v_5, v_6, v_7, v_8 \rangle \subset \\ \subset \langle v_1, v_2, v_4, v_5, v_6, v_7, v_8 \rangle \subset V_9/V_1.
\end{multline}
given by the desingularization of $Y_6$. As it turns out, any flag in the rational normal curve
\[
V_3/V_1 = \langle v_8, \alpha v_4 + \beta v_7 \rangle \subset
V_6/V_1 = \langle v_4, v_5, v_7, v_8, \alpha v_2 + \beta v_6 \rangle
\]
is contained in $p^{-1}(V_1)$ since $\omega \in V_{199}+V_{369}+V_{666}$, hence the conclusion follows if we can prove that it is contained in $E$. On any such flag, \eqref{uniqueflag} induces flags $A_1\subset A_2$, $B_1\subset B_2
\subset B_3$ and  $C_1\subset C_2\subset C_3$ such that 
\[\bar\omega \in A_1\otimes B_2\otimes C_3+A_2\otimes B_1\otimes C_2+(A_1\otimes B_3+A_2\otimes B_2)\otimes C_1.\]
Then it is easy to see that the length three subscheme we get in  $\PP(C_3^\vee)$ has multiplicity 
two at the point defined by the hyperplane $C_2$ (note that this point is exactly the hyperplane $V_8$, 
uniquely defined by $V_1\in D_{Y_6}$). Since $E$ is precisely the locus of non-reduced 
schemes, we are done. 
\end{proof}
\end{theorem}

Let $P\in A$ be a hyperplane in $\PP(V_9)$, and let $P_4$ be the four-dimensional subspace of $P$ such that $\omega(P,\cdot,\cdot)\in \wedge^2 P_4$. Then:

\begin{proposition}
\label{coverbyp3}
$D_{Y_3}$ is covered by a family of $\PP^3$ parametrized by $A$. More precisely, for any point $P\in A$, we have that $\PP(P_4)\subset D_{Y_3}$.

\begin{proof}
As $P\in A$ and by the definition of $P_4$, we know that 
\[\omega\in \wedge^3 P + V_9\wedge (\wedge^2 P_4).\]
In order to show that $\PP(P_4)\subset D_{Y_3}$, we need to prove that for any $V_1\subset P_4$, there exist $U_2\subset U_5\subset V_9$ such that $V_1\subset U_2$ and $\omega\in U_2\wedge (\wedge^2 V_9) +V_9\wedge (\wedge^2 U_5)$. Indeed, if this happens, then $\omega$ modulo $V_1$ belongs to the total space of the vector bundle which gives a Kempf collapsing of $Y_3$ inside $V_9/V_1$ (see Table \ref{kempfTable}), and therefore $V_1\in D_{Y_3}$.

Let $V_1\subset P_4$. We construct $U_5$ as a subspace of $P$. Moreover, $\omega(P,\cdot,\cdot)$ is a two-form on $P_4$, and therefore we can consider the orthogonal $V_1^{\perp_\omega}\subset P_4$ of $V_1$ inside $P_4$ with respect to this two-form. We construct $U_2$ as a subspace of $V_1^{\perp_\omega}$ containing $V_1$. The parameter space for $U_2$ is then $\PP(V_1^{\perp_\omega}/V_1)$ and the parameter space for the pair $U_2\subset U_5$ is the Grassmannian bundle $\Gr(3,P/U_2)$ over $\PP(V_1^{\perp_\omega}/V_1)$, a variety of dimension ten.

Asking that $V_1\subset U_2\subset V_1^{\perp_\omega}$ implies that $\omega(P,\cdot,\cdot)\in U_2\wedge P_4$. Therefore we have that $\omega\in \wedge^3 P + U_2\wedge (\wedge^2 V_9)$. Let us consider the element $\tilde{\omega}\in \wedge^3(P/U_2)$ induced by $\omega$. Then $\omega\in U_2\wedge (\wedge^2 V_9) +V_9\wedge (\wedge^2 U_5)$ if and only if $\tilde{\omega}\in \wedge^2 (U_5/U_2) \wedge P/U_2$. Over our parameter space, $U_5/U_2$ is parametrized by the rank three tautological bundle $\cU$ over $\Gr(3,P/U_2)$. As a consequence, requiring that $\tilde{\omega}\in \wedge^2 \cU\wedge (P/U_2)$ is the same as asking that the induced section $\bar{\omega}$ of the vector bundle $F:=\wedge^3 (P/U_2)/(\wedge^2\cU\wedge (P/U_2))$ vanishes. $F$ is a rank ten vector bundle over the ten-dimensional parameter space, and the zero locus of its general section $\bar{\omega}$ parametrizes the pairs $U_2, U_5$ such that $\omega\in U_2\wedge (\wedge^2 V_9) +V_9\wedge (\wedge^2 U_5)$. This zero locus consists in general of
\[
c_{10}(F)=2
\]
points, as a computation with \cite{Macaulay2} shows; as it is nonempty, there exist $U_2\subset U_5$ with the required properties, and $V_1\in D_{Y_3}$. This concludes the proof. (Note that 
the existence of different flags comes from the fact that the Kempf collapsing we used has degree $d>1$, see Table \ref{kempfTable}: the above computation actually shows that $d=2$, as stated in Remark \ref{nonbiry3}.)
\end{proof}

\end{proposition}

\begin{remark}
The singular locus of $SU_C(3)$ can be identified with $D_{Y_3}$. This singular locus 
is known to coincide with the set of strictly semistable rank three vector bundles with trivial determinant. A generic point of this set is a bundle $L\oplus E$, where $L$ is a line bundle (or a point of $A$) and $E$ is a rank two vector bundle such that $\det(E)=L^{-1}$. Therefore, having fixed $L$, this set contains \[
\{ E \mbox { of rank $2$ s.t.\ } \det(E)=L^{-1}\} \cong SU_C(2),
\]
which is a $\PP^3$, see Section \ref{modvecbun}. This gives the family of $\PP^3$'s parametrized by $A$ covering $D_{Y_3}$ exhibited in Proposition \ref{coverbyp3}. Moreover, each $\PP^3$ contains a copy of $\KU$, the Kummer surface associated to $A$: this gives the family of $\KU$ parametrized by $A$ covering $D_{Y_4}\cong \Sigma$.
\end{remark}

\begin{corollary}
$D_{Y_6}$ is not normal, and $A$ is its normalization.
\begin{proof}

We know that $D_{Y_6}$ is the singular locus of $D_{Y_4}$, so by Theorem \ref{sigma} it coincides with the set of triples of the form $(P,P,-2P)$ in $\Sigma$. In particular there is a bijective morphism $A\ra D_{Y_6}$, which implies that the normalization of $D_{Y_6}$ is isomorphic to $A$.

There just remains to prove that the singular locus of $\Sigma$ is not normal. Let us consider the following commutative diagram
\[
\xymatrix{
\Sigma \times A \ar[r]^-{\alpha} & A^{(3)}\\
\alpha^{-1}(\Sing (A^{(3)})) \ar[r] \ar@{^{(}->}[u] & \Sing (A^{(3)}) \ar@{^{(}->}[u]
}
\]
where $\alpha:([P,Q,R],t)\mapsto[P+t,Q+t,R+t]$. The preimage of a point is
\[
\alpha^{-1}([P',Q',R'])=\{([P'-s,Q'-s,R'-s],s) \mbox{ s.t.\ } 3s=P'+Q'+R'\};
\]
in particular, $\alpha$ is a $81:1$ \'etale cover of $A^{(3)}$, and induces an \'etale cover of $\Sing (A^{(3)})$ by the singular locus of $\Sigma \times A$. Consider the following diagram:
\[
\xymatrix{
\Delta_{12}\cup\Delta_{13}\cup\Delta_{23}  \ar@{->>}[d] \ar@{^{(}->}[r] & A^3 \ar@{->>}[d]\\
\Sing (A^{(3)}) \ar@{^{(}->}[r] & A^{(3)}
}
\]
where, e.g., $\Delta_{12}=\{(P,P,Q) \in A^3\}$ is the first diagonal. The map on the left is generically $3:1$, while the map on the right is generically $6:1$. The restriction $\Delta_{12}\rightarrow \Sing (A^{(3)})$ is a birational finite morphism, hence it is an isomorphism if $\Sing (A^{(3)})$ is normal.

Locally, we have $A^3 \cong (\mathbb{C}^2)^3$ and $\cO_{A^3} \cong R:=\mathbb{C}[x_1,y_1,x_2,y_2,x_3,y_3]$. The above diagram induces the following commutative diagram
\[
\xymatrix{
\cO_{\Delta_{12}} & \ar@{->>}[l] R \\
\cO_{\Sing (A^{(3)})} \ar@{^{(}->}[u]^-{\beta} & \ar@{->>}[l] R^{S_3} \ar@{^{(}->}[u]
}
\]
To conclude, it is enough to show that $\beta$ is not an isomorphism. Locally, $\cO_{\Delta_{12}} \cong R/(x_1-x_2,y_1-y_2)$, while by construction $\cO_{\Sing (A^{(3)})}$ is the quotient of $R^{S_3}$ by the homogeneous ideal $(x_1-x_2,y_1-y_2)^{S_3}$. Since for instance $x_3 \in \cO_{\Delta_{12}}$ but it is not in the image of $\beta$, the conclusion follows.
\end{proof}
\end{corollary}

\begin{corollary}
	The orbit closure $Y_4 \subset \wedge^3 V_8$ is singular along $Y_6$. The orbit closure $Y_6 \subset \wedge^3 V_8$ is non-normal along $Y_8$.
\end{corollary}
We observe that the last statement agrees with and specifies the claim in \cite{weymanE8} about the non-normality of $Y_6$.

\begin{remark}
\label{remrestriction}
One can observe that the isomorphism $D_{Y_4} \ra \Sigma$ constructed in Theorem \ref{sigma} restricts to the birational map $D_{Y_6}\dashrightarrow A$ described in Proposition \ref{propbirA}. A point $x$ of $D_{Y_6}\setminus D_{Y_8}$ corresponds to $[P,P,Q]$, where $P \in \PP(V_9^\vee)$ is the hyperplane defined by the preimage of $x$ in the desingularization of $D_{Y_6}$. Similarly, a point $x' \in D_{Y_8}$ corresponds to $[P',P',P']$, where $P' \in \PP(V_9^\vee)$ is the hyperplane defined by the preimage of $x'$ in the desingularization of $D_{Y_8}$.
\end{remark}

\subsection{The nested Kummer fourfold and the Hilbert scheme}
Let us consider now the nested Hilbert scheme $\Hilb^{2,3}(A)$ parametrizing pairs
$(S,S')$, where $S$ is a length two subscheme of $A$, and $S'\supset S$ a length three 
subscheme. Such a nested Hilbert scheme is known to be smooth. Moreover,
it admits an action of $A$ by translation, compatible with the sum map. So all the 
fibers of the latter are equivalent, and smooth. We denote them by $\KN$,
the nested Kummer fourfold. By restriction from the Hilbert schemes, we get 
a triple cover $\KN\ra \KK$, branched over the exceptional divisor, 
and also a morphism $\KN\ra A$ defined by taking the residual scheme. 

In our situation, $T$ is a triple cover of $D_{Y_4}$, which is birational to 
$Z\simeq \KK$. So the fiber product of $Z$ with $T$ over $D_{Y_4}$ will be
a triple cover of $\KK$, and we can expect it to be isomorphic to 
$\KN$. Rather than taking formally the direct product, we define $X\subset 
F(1,3,5,6,7,V_9)$ as parametrizing the flags $V_1\subset V_3\subset V_5\subset V_6\subset V_7\subset V_9$ such that $\omega$ belongs to $V_{199}+V_{359}+V_{377}+V_{666}$. Just like 
$T$ and $Z$, for $\omega$ generic this is a smooth fourfold. Since $V_{199}+V_{359}+
V_{377}+V_{666}$ is exactly the intersection of $V_{199}+V_{369}+V_{666}$ and $V_{199}+V_{559}+V_{777}$, $X$ admits natural projections to $Z$ and $T$:
\begin{equation}
\label{Adiag1}
\xymatrix @C=2pc @R=0.8pc{
& & \rule{1pt}{0pt}X \ar[dl]_{{3:1}} \ar[dr]^{\mbox{\small birat}}  \\
& Z \ar[dr]_{\mbox{\small desing}} &  & T \ar[dl]^{{3:1}}  \\
 & & \Sigma & 
} 
\end{equation}
Beware that the degree three morphisms in this diagram are only generically finite. In view of Theorem \ref{isoK23}, we can give a precise statement concerning $X \ra Z$:

\begin{lemma}\label{almostfinite}
The positive dimensional fibers of $X\ra Z$ are $81$ projective lines.

\begin{proof}
We will prove that the positive dimensional fibers are $81$ projective lines which are contracted to $81$ points in $Z$, whose image via the desingularization $Z \ra \Sigma$ is precisely $D_{Y_8}$. The reason behind this phenomenon is further clarified in Remark \ref{rmpositivefibers} below.

A point $z$ of $Z$ is a flag $V_1\subset V_3\subset V_6\subset V_9$ such that $\omega$ belongs to
$M=V_{199}+V_{369}+V_{666}$. A point $p$ of $X$ above $z$ is a pair of subspaces $(V_5,V_7)$, with 
$V_1\subset V_3\subset V_5\subset V_6\subset V_7\subset V_9$, such that $\omega$ belongs to 
$N=V_{199}+V_{359}+V_{377}+V_{666}$. This is a subspace of $M$, and 
\[M/N\simeq V_{369}/(V_{359}+V_{169}+V_{367})=V_3/V_1\otimes V_6/V_5\otimes V_9/V_7.\] 
Let $A_2=V_3/V_1$. 
The pairs $(V_5,V_7)$ are parametrized by the product of Grassmannians $\Gr(2,V_6/V_3) \times \Gr(1,V_9/V_6)\simeq \PP^2\times \PP^2$.
On this variety, $\omega$ defines a global section of the rank four vector bundle $A_2\otimes V_6/V_5\otimes V_9/V_7
=A_2\otimes \cQ_1\otimes \cQ_2$, and the fiber over $z$ identifies with the zero locus of this section. Here 
we have denoted by $\cQ_1$ and $\cQ_2$ the quotient bundles, of rank one and two, over $\Gr(2,V_6/V_3)$ and $\Gr(1,V_9/V_6)$. 
An easy computation shows that $c_4(A_2\otimes \cQ_1\otimes \cQ_2)=c_2(\cQ_1\otimes \cQ_2)^2=3$, confirming that 
the general fiber consists in three points. 

There remains to identify the infinite fibers. For this we need to analyze when a section of 
$A_2\otimes \cQ_1\otimes \cQ_2$ on $\Gr(2,B_3)\times \Gr(1,C_3)$ vanishes in positive dimension, where $B_3=V_6/V_3$ and $C_3=V_9/V_6$. 
Such a 
global section is an element $\gamma$ of $A_2\otimes B_3\otimes C_3$; as $A_2\otimes B_3\otimes C_3=M/V_{199}+V_{339}+V_{666}$, we have
\[
\gamma = \omega \, (\modulo V_{199}+V_{339}+V_{666}).
\]
We consider $\gamma$ as a 
family $\Gamma$ of linear maps from $B_3^\vee$ to $C_3$. This section vanishes at $(U_2,U_1)$, 
where $U_2\subset B_3$ and  $U_1\subset C_3$, if and only if $\gamma$ belongs to $A_2\otimes (U_2\otimes C_3+B_3\otimes U_1)$. In other words, all the linear maps in $\Gamma$ send the 
line $U_2^\perp$ to the line $U_1$.

The classification of pencils of $3\times 3$-matrices is well-known: there are exactly seventeen orbits (see e.g.\ \cite[5.4]{weymanE6}). A straightforward check shows that the maximal orbits such that $\gamma$ vanishes in infinitely many pairs $(U_2,U_1)$ are those named as $\cO_{14}$, $\cO_{13}$, $\cO_{11}$, $\cO_{10}$, whose elements can be written respectively as follows:
\begin{enumerate}[label=\roman*.]
\item $a_1\otimes b_1\otimes c_1+
a_2\otimes b_2\otimes c_2+a_1\otimes b_3\otimes c_3$ for some $a_i\in A_2$, $b_j\in B_3$,
$c_k\in C_3$;
\item $a_1\otimes b_1\otimes c_1+a_2\otimes b_2\otimes c_2+a_1\otimes b_2\otimes c_3+
a_2\otimes b_3\otimes c_1$ for some $a_i\in A_2$, $b_j\in B_3$,
$c_k\in C_3$.
\item $a_1\otimes b_1\otimes c_1+a_2\otimes b_2\otimes c_2+a_1\otimes b_2\otimes c_3+
a_2\otimes b_1\otimes c_3$ for some $a_i\in A_2$, $b_j\in B_3$,
$c_k\in C_3$;
\item $a_1\otimes b_1\otimes c_1+a_2\otimes b_2\otimes c_2+a_1\otimes b_3\otimes c_2+
a_2\otimes b_3\otimes c_1$ for some $a_i\in A_2$, $b_j\in B_3$,
$c_k\in C_3$.
\end{enumerate}
We will show that $\omega$ being generic, $\gamma$ will never be of any of the first three types. By contradiction, we will prove that $\gamma$ being of those special types would force the class of $\omega$ modulo $V_1$ to belong to some orbit 
closure $Y$ 
of codimension bigger than eight in $\wedge^3(V_9/V_1)$. In other words, we would 
get a point in an ODL $D_Y\subset\PP (V_9)$, which we know to be empty for 
a generic $\omega$.
\begin{enumerate}[label=\roman*.]
\item In case {i.},
we can find 
$A_1\subset A_2$, $B_1\subset B_3$, $C_1\subset C_3$ such that  $\gamma$
belongs to $A_1\otimes B_3\otimes C_3+A_2\otimes B_1\otimes C_1$. This means that 
we can find $V_2$, $V_4$, $V_7$, with $V_1\subset V_2\subset V_3\subset V_4
\subset V_6\subset V_7\subset V_9$ such that $\omega$ belongs to $V_{199}+V_{269}+V_{347}+V_{666}$. Modding out by $V_1$ and letting $U_i=V_{i+1}/V_1$, 
we get a point in the total space of the vector bundle $U_{158}+U_{236}+U_{555}$ 
over $F(1,2,3,5,6,U_8)$.
Note that this vector bundle is a subbundle of $U_{168}+U_{666}$, which has rank $30$ over the $17$-dimensional flag manifold $F(1,6,U_8)$.
So it will collapse to an orbit closure of dimension at most $30+17<48$ in 
$\wedge^3U_8$. 
\item In case {ii.}, we observe that we can find $B_1\subset B_3$ and $C_1\subset C_3$
such that $\gamma$ belongs to $A_2\otimes B_1\otimes C_3+A_2\otimes B_3\otimes C_1$.
This means that we can find $V_4$, $V_7$, with $V_3\subset V_4
\subset V_6\subset V_7\subset V_9$ such that $\omega$ belongs to $V_{199}+V_{349}+V_{367}+V_{666}$. Modding out by $V_1$, 
we get a point in the total space of the vector bundle $U_{238}+U_{256}+U_{555}$ 
over $F(2,3,5,6,U_8)$. The latter flag manifold has dimension $25$, and the vector bundle
has rank $23$. But note that $U_{238}+U_{256}+U_{555}\subset U_{338}+U_{666}$, and that
the vector bundle $U_{338}+U_{666}$ has rank $26$ over the $21$-dimensional flag manifold
$F(3,6,U_8)$. So again it will collapse to an orbit closure of dimension at most $26+21<48$ in 
$\wedge^3U_8$.
\item In case {iii.},
there exists $B_2\subset B_3$ such that $\gamma$ belongs to $A_2\otimes B_2\otimes C_3$.
So there exists $V_5$, with $V_3\subset V_5\subset V_6$, such that $\omega$ belongs to 
$V_{199}+V_{359}+V_{666}$.  Modding out by $V_1$ as before,
we get a point in the total space of the vector bundle $U_{248}+U_{555}$
over $F(2,4,5,U_8)$. The latter flag manifold has dimension $23$, and the vector bundle
has rank $25$, so it would seem to collapse to a codimension $8$ orbit closure. 
But notice that 
$U_{248}+U_{555}\subset U_{448}+U_{555}$, and that now $U_{448}+U_{555}$ is a rank $28$ 
vector bundle over $F(4,5,U_8)$, whose dimension is $19$. So the collapsing will have for
image an orbit closure of dimension at most $28+19<48$.
\end{enumerate}

So we only remain with case {iv.}. Observe that it occurs if and only if there exists $C_2\subset C_3$ such that $\gamma$ belongs to $A_2\otimes B_3\otimes C_2$. So there exists $V_8$, with $V_6\subset V_8\subset V_9$, such that $\omega$ belongs to $V_{199}+V_{339}+V_{368}+V_{666}$, and its class modulo $V_1$ is contained in some $U_{228}+V_{257}+U_{555}$. This is the vector bundle 
on $F(2,5,7,U_8)$ that desingularizes $Y_8$. In particular $V_1$ defines one of the 
$81$ points of $D_{Y_8}\subset\PP(V_9)$, the flag $V_3\subset V_6\subset V_8$ 
is uniquely defined by $V_1$ and there is a uniquely defined lift of 
$V_1$ in $Z$. An easy computation shows that the fiber in $X$ of this lift is a projective line, as claimed, and no further degeneration of $\omega$ can occur for $\omega$ generic.
\end{proof}
\end{lemma}

\begin{theorem}\label{isoK23}
$X$ is isomorphic to the nested Kummer fourfold $\KN$.

\begin{proof} We would like to lift the isomorphism between $Z$ and $\KK$:
\begin{equation}
\label{K23diag}
\xymatrix @C=2pc @R=0.8pc{
X\ar[d]\ar@{.>}[r]^-? &  \KN \ar[d]\\
Z\ar[r]^-\simeq & \KK
} 
\end{equation}
By definition, $X$ parametrizes the flags $V_1\subset V_3\subset V_5\subset V_6\subset V_7\subset V_9$ such that 
\[\omega\in V_{199}+V_{359}+V_{377}+V_{666}.\]
Consider a hyperplane $V_8=\Ker(\phi)$ containing $V_7$. It defines a point in $A$ 
if the contraction of $\omega$ by $\phi$ has rank four. Let $v_8, v_9$ be vectors in $V_9$, 
independent modulo $V_7$. Modulo $V_{777}$, which is killed
by $\phi$, we can write $\omega = v_1\wedge \alpha +\beta_8 \wedge v_8 +\beta_9\wedge v_9$,
for some $\alpha\in\wedge^2V_9$ and $\beta_8, \beta_9\in V_3\wedge V_5$. Therefore 
\[\phi\lrcorner \omega = -v_1\wedge (\phi\lrcorner \alpha)+\phi(v_8)\beta_8+\phi(v_9)\beta_9.\]
Consider the pencil $\langle\beta_8, \beta_9\rangle$. If we mod out by $V_1\wedge V_5$, 
we get $\langle\bar\beta_8, \bar\beta_9\rangle$ in $\wedge^2(V_5/V_1)$, hence in general
a pencil that cuts the Grassmannian of rank two tensors along a length two subscheme. 
Substracting $v_1\wedge (\phi\lrcorner \alpha)$ to a rank two tensor yields a tensor of
rank at most four. So we get a rational map $X\ra \Hilb^2 (A)$. 

What could prevent it to be regular? First, it could happen that when we mod out by 
$V_1\wedge V_5$, the pencil $\langle\beta_8, \beta_9\rangle$ collapses. In other 
words, $\beta_8$ and $\beta_9$ could be proportional up to $V_1\wedge V_5$. Then 
we may suppose that $v_9=0$. But this would mean that modulo $V_1$, $\omega$ depends
only on seven variables, a condition that inside $\wedge^3\CC^8$ defines an orbit of 
codimension $14>8$. So this cannot happen.

Second, the projected pencil $\langle\bar\beta_8, \bar\beta_9\rangle$ could be
contained in the Grassmannian of rank two tensors. But then the pencil of hyperplanes 
that contain $V_7$ would be contained in $A$. Since an abelian surface cannot contain
any line, this cannot happen either. 

Combining the regular map $X\ra \Hilb^2( A)$ with the projection $X\ra Z\simeq \KK$, we get a morphism $X\ra \KK\times \Hilb^2 (A)$ whose image is by construction contained in, hence equal to,
$\KN$. Since the projections from $X$ to $Z$ and from $\KN$
to $\KK$ are both generically finite of degree three, we get a birational morphism 
from $X$ to $\KN$. But by Lemma \ref{almostfinite} the exceptional locus of 
this birational morphism is at most one-dimensional. So it has to be an isomorphism.
 \end{proof}
\end{theorem}

\begin{remark}
\label{rmpositivefibers} The positive dimensional fibers of the projection map
from $\KN$ to $\KK$ live above the $81$ three-torsion points of $A$. 
Indeed, if $P$ is such a point, then the fat point defined by $I(P)^2$ is a length three
subscheme that defines a point in $\KK$. Since this subscheme contains all the 
length two schemes supported at $P$, we get a fiber isomorphic to $\PP(T_PA)\simeq\PP^1$. 
All the other length three schemes supported at $P$ are curvilinear, hence define a unique
tangent. In particular we get an identification of $D_{Y_8}$ with $A[3]$, provided we have fixed an origin (in the preimage of  $D_{Y_8}$ in $A$ via the map $A \ra D_{Y_6}$).
\end{remark}

\begin{theorem}\label{isoThilbert}
$T$ is isomorphic to $\Hilb^2 (A)$. 

\begin{proof} In the proof of Theorem \ref{isoK23} we constructed a morphism $\eta_X: X\ra \Hilb^2 (A)$. 
In fact the construction shows that this is the composition of a morphism from $\eta_T : T
\ra \Hilb^2 (A)$ with the projection $X\ra T$. Since this morphism is birational, as well 
as the projection $\KN\ra \Hilb^2 (A)$, $\eta_T$ is birational. Since $T$ has trivial
canonical bundle, such a birational morphism must be an isomorphism.\end{proof}
\end{theorem}

\begin{remark}
The group structure of $A$ allows us to define a surjective morphism $A\times \KU$
to $\Hilb^2 (A)$ which is an \'etale cover of degree sixteen. This is the \'etale cover
whose existence is predicted by the Beauville--Bogomolov decomposition. 
\end{remark}

\section{On the group structure of \texorpdfstring{$A$}{A}}

In this section we geometrically describe the group structure of $A$. This is an analogue of the usual description of the group structure of a 
plane cubic curve from its intersection with lines. It is worth mentioning that, in a 
different context, Donagi \cite{DonagiGroup} provided a geometric characterization of the group law for the $n$-dimensional abelian variety parametrizing the $(n-1)$-dimensional linear subspaces of the intersection of two general quadrics in $\PP^{2n+1}$, which is known to be the jacobian of a hyperelliptic curve of genus $n$.
\smallskip

Recall what we have established so far. If we choose two distinct points $P,Q$ of 
$A\subset \PP(V_9^\vee)$, the corresponding point $z$ in $\Hilb^2 (A)\simeq T$ maps to 
a point $[V_1]\in D_{Y_4}\subset \PP(V_9)$. If this point is not on $D_{Y_6}$, 
it defines a flag $V_1\subset V_3\subset V_6\subset V_9$ such that 
$\omega$ belongs to $V_{199}+V_{369}+V_{666}$. Moreover, its three preimages 
$z, z', z''$ in $T$ yield additional subspaces $(V_5, V_7)$, $(V'_5, V'_7)$, 
$(V''_5, V''_7)$, with 
\[V_1\subset V_3\subset V_5, V'_5, V''_5\subset V_6\subset V_7, V'_7, V''_7\subset V_9,\]
such that $\omega$ belongs to $V_{199}+V_{559}+V_{777}$, and to the corresponding 
spaces with $(V_5, V_7)$ replaced by $(V'_5, V'_7)$ and $(V''_5, V''_7)$. 
Since $V_7=P\cap Q$, this implies that if we contract $\omega$ by an equation of
the hyperplane $P$ and an equation of the hyperplane $Q$, we get a vector in $V_1$. 
With a slight abuse of notation, we write 
\[[\omega(P, Q,\cdot )]=[V_1]\in D_{Y_4}\subset \PP(V_9).\]
This yields a simple description of the map from $\Hilb^2 (A)\simeq T$ to $D_{Y_4}$. 
Moreover, this is enough to characterize the point $R\in A$ such that $(P,Q,R)$ belongs
to $\Sigma$:

\begin{proposition}
\label{lemmaPQR}
Let $P, Q\in \PP(V_9^\vee)$ be general points of $A$. Then the unique point $R\in A$ 
such that $(P,Q,R)$ belongs to $\Sigma$ is characterized by the condition 
\[
[\omega(P, R,\cdot )]= [\omega(Q, R,\cdot )]=[\omega(P, Q,\cdot )].
\]

\begin{proof}
The previous remarks show that $R$ verifies the required condition. There remains to 
prove that it is uniquely characterized by it. 

Let $V_7:=P\cap Q $, and let us choose $v_P\in Q\setminus P$ and $v_Q\in P \setminus Q$. Let us decompose 
$\omega$ with respect to the direct sum $V_9=\CC v_P\oplus \CC v_Q\oplus V_7$, as 
\[
\omega = v_P\wedge v_Q\wedge v_1 + v_P\wedge \alpha +v_Q\wedge\beta+\sigma,
\] 
with $\alpha, \beta\in \wedge^2 V_7$ and $\sigma\in \wedge^3 V_7$.
In particular $v_1$ generates $V_1$.
Since $P$ belongs to $A$, 
$v_Q\wedge v_1 + \alpha$ has rank (at most) four (and since also $Q$ belongs to $A$, 
$v_Q\wedge v_1 - \beta$ also has rank (at most) four). This means that $\alpha$ itself
has rank at most four, and $v_Q\wedge v_1\wedge \alpha\wedge\alpha=0$, or equivalently
$v_1\wedge \alpha\wedge\alpha=0$.

\begin{lemma}
$v_1\wedge \alpha\wedge\alpha=0$ if and only if there exist $u,v,w$ such that $\alpha=
v_1\wedge u + v\wedge w$. 

\begin{proof} If $\alpha$ has rank six or more, then $v_1\wedge \alpha\wedge\alpha=0$ would imply
$v_1=0$, which is not the case. If $\alpha$ has rank two, $\alpha\wedge\alpha=0$. So suppose
that $\alpha$ has rank exactly four, which means that there exists a unique four-plane 
$L$ such that $\alpha$ belongs to $\wedge^2L$. Then $\alpha\wedge\alpha$ is a generator of
$\wedge^4L$, and $v_1\wedge \alpha\wedge\alpha=0$ means that $v_1$ belongs to $L$. The 
conclusion easily follows, since if we choose a generic vector $u$ in $L$, the line
generated by $\alpha$ and $v_1\wedge u$ in $\PP (\wedge^2L)$ will meet the quadric 
of rank two tensors at another point.
\end{proof}
\end{lemma}

Applying this Lemma also to $v_Q\wedge v_1 - \beta$, we deduce that there exist $u,v,w , u', v', w'\in V_7$ such that 
\begin{equation*}\label{PQ}
\omega=v_P\wedge v_Q\wedge v_1 + v_P\wedge (v_1\wedge u + v\wedge w) +v_Q\wedge(v_1\wedge u' +v'\wedge w') +\sigma.
\end{equation*}
Generically, $v_P, v_Q, v_1, u,v,w , u', v', w'$ is a basis of $V_9$.  
Note that $\omega(P, R,\cdot )$ is the contraction of $v_Q\wedge v_1 + v_1\wedge u + v\wedge w$ by $R$ (considered as a linear form). In particular it will be proportional to $v_1$ if and only if $v_1, v, w$ belong to the hyperplane $R$. Similarly  $\omega(Q, R,\cdot )$ is generated by $v_1$ 
if and only if $v_1, v', w'$ belong to $R$. So $R\supset U_5=\langle v,w,v',w',v_1 \rangle$.
If we let $U_4=\langle v_P, v_Q, u, u' \rangle$, we are thus looking for  $R\in \PP(U_4^\vee)
\cong \PP^3$ such that $\omega(R,\cdot,\cdot)$ has rank four. 

Let us decompose $\sigma$ further with respect to the decomposition $V_7=\CC u\oplus \CC u'
\oplus U_5$: there exist $\sigma_{0}\in U_5$, $\sigma_u, \sigma_{u'}\in \wedge^2U_5$ and $\tilde{\sigma}\in \wedge^3U_5$ such that
\[
\sigma=u\wedge u'\wedge \sigma_0 + u\wedge \sigma_u +u'\wedge \sigma_{u'}+\tilde{\sigma}.
\]
As a consequence, we get 
\[
\omega(R,\cdot,\cdot)=  a\wedge v_1 + b\wedge\sigma_0+\tau,
\]
where $a=R(v_P+u')v_Q + R(u-v_Q)v_P-R(v_P)u-R(v_Q)u'$ and $b=R(u)u'-R(u')u$ belong
to $U_4$, while  
\[
\tau=R(v_P)v\wedge w+R(v_Q)v'\wedge w' + R(u)\sigma_u+R(u')\sigma_{u'}
\]
belongs to $\wedge^2U_5$ (here again we denoted by the same letter $R$ a linear form 
whose kernel is the hyperplane $R$). 
If $a$ and $b$ are dependent and for example $b=0$, we need 
that $R(u)=R(u')=0$. Then $\omega(R,\cdot,\cdot)=  a\wedge v_1+R(v_P)v\wedge w
+R(v_Q)v'\wedge w'$ never has rank four or less, unless it is zero. 
The case where $b\ne 0$ is similar. 

If $a$ and $b$ are independent, since generically $v_1$ and $\sigma_0$ are independent, 
the only way for $\omega(R,\cdot,\cdot)$ to have rank at most four is that 
$v_1\wedge \sigma_0\wedge \tau=0$. Since the map $v_1\wedge \sigma_0\wedge : \wedge^2V_5\longrightarrow
\wedge^4V_5$ has rank three, this yields three linear conditions that determine $(R(v_P), 
R(v_Q), R(u), R(u'))$ uniquely up to scalar. So the hyperplane $R$ is uniquely determined.
\end{proof}
\end{proposition}

\medskip
The point $[\omega(P, Q,\cdot )]=[V_1]\in D_{Y_4}$ should really be thought of as the line joining $P$ and $Q$, in analogy 
to the line joining two points on a plane cubic, and that defines a unique third 
point. From this perspective, the space of ``lines'' is $D_{Y_4}$, or rather its 
desingularization $Z=\KK$. 

Once we have chosen an origin $O$ of $A$, exactly as for plane cubics we can then recover 
the group structure on $A$
by applying  Proposition \ref{lemmaPQR} twice: starting from two general points $P,Q\in A$, 
we first find the point $R$ such that $(P,Q,R)$ belongs to $\Sigma$; then from the two
points $O,R\in A$, we deduce the sum $P+Q\in A$.

Finally, the isomorphism $A \to A^\vee$ described in Section \ref{altTriv}, which associates to each point $x \in A$ the curve $C_x$ obtained from the intersection of $A$ with $\mathbf{P}(\ker(\omega(x,\cdot,\cdot)))$, allows us to reinterpret Proposition \ref{lemmaPQR} as follows. 

\begin{proposition}
	\label{lemmaPQRcurves}
	Let $P, Q \in \PP(V_9^\vee)$ be general points of $A$. Then $(P,Q,R)$ belongs to $\Sigma$ for $R \in A$ if and only if there exists a hyperplane $H \subset \PP(V_9^\vee)$ cutting out $C_P \cup C_Q \cup C_R$ on $A$.
	\begin{proof}
		Let us consider $[\omega(P, Q,\cdot )] \in \PP(V_9)$. The corresponding hyperplane in $\PP(V_9^\vee)$ contains, by construction, both $\mathbf{P}(\ker(\omega(P,\cdot,\cdot)))$ and $\mathbf{P}(\ker(\omega(Q,\cdot,\cdot)))$ (which meet in the two points $C_P \cap C_Q$ and indeed generate a $\PP^7$), hence it cuts out on $A$ a reducible curve containing $C_P \cup C_Q$. If $(P,Q,R)$ belongs to $\Sigma$, by Proposition \ref{lemmaPQR} the same applies for the pairs $(P,R)$ and $(Q,R)$, hence the residual curve is $C_R$. Conversely, suppose that $H \cap A=C_P \cup C_Q \cup C_R$. Then $H$ necessarily coincides with the hyperplane corresponding to $[\omega(P, Q,\cdot )] \in \PP(V_9)$. By permuting the points, the conclusion follows by Proposition \ref{lemmaPQR}.
	\end{proof}
\end{proposition}
A posteriori, $D_{Y_4}$ parametrizes all hyperplanes cutting $A$ in the union of three degree six curves, possibly counted with multiplicities.

\makeatletter
\providecommand\@dotsep{5}
\makeatother

\listoftodos 

\end{document}